\documentclass{amsart}
\usepackage{graphicx}
\usepackage{amsmath}
\usepackage{amssymb}
\usepackage{xcolor}
\usepackage{mathtools}
\usepackage{geometry}
\newtheorem{theorem}{Theorem}[section]
\newtheorem{proposition}[theorem]{Proposition}
\newtheorem{Def}[theorem]{Definition}
\newtheorem{Cor}[theorem]{Corollary}
\newtheorem{Lemma}[theorem]{Lemma}
\newtheorem{rem}[theorem]{Remark}

\title[Euler factors of equivariant $L$--functions of Drinfeld modules]{Euler factors of equivariant $L$--functions\\ of Drinfeld modules and beyond}
\author{Cristian D. Popescu}
\author{Nandagopal Ramachandran}
\address{Dept. of Mathematics, University of California at San Diego, San Diego, CA 92093, USA}
\email{cpopescu@ucsd.edu}
\email{naramach@ucsd.edu}



\date{}

\begin{document}
\begin{abstract} In \cite{FGHP}, the first author and his collaborators proved an equivariant Tamagawa number formula for the special value at $s=0$ of a Goss--type $L$--function, equivariant with respect to a Galois group $G$, and associated to a Drinfeld module defined on $\Bbb F_q[t]$ and over a finite, integral extension of $\Bbb F_q[t]$. The formula in question was proved provided that the values at $0$ of the Euler factors of the equivariant $L$--function in question satisfy certain identities involving Fitting ideals of certain $G$--cohomologically trivial, finite $\Bbb F_q[t][G]$--modules associated to the Drinfeld module. In \cite{FGHP}, we prove these identities in the particular case of the Carlitz module. In this paper, we develop general techniques and prove the identities in question for arbitrary Drinfeld modules. Further, we indicate how these techniques can be extended to the more general case of higher dimensional abelian $t$--modules, which is relevant in the context of the proof of the equivariant Tamagawa number formula for abelian $t$--modules given by N. Green and the first author in \cite{Green-Popescu}. This paper is based on a lecture given by the first author at ICMAT Madrid in May 2023 and builds upon results obtained by the second author in his  PhD thesis \cite{Ramachandran-thesis}.

\end{abstract}

\maketitle

\tableofcontents
\section{Introduction} Let $E/F/k$ be a tower of finite field extensions, where $k:=\Bbb F_q(t)$, $q$ is a power of a prime $p$, and $E/F$ is Galois, of abelian Galois group $G$. Further, let $\mathcal O_F$ and $\mathcal O_K$ denote the integral closure of $A:=\Bbb F_q[t]$ in $F$ and $K$, respectively, and consider a Drinfeld module $E$ defined on $A:=\Bbb F_q[t]$ and with coefficients in $\mathcal O_F$. 

In \cite{FGHP}, the authors associated to the set of data $(K/F/k, E)$ a $G$--equivariant Goss--type $L$--function 
$$\Theta^E_{K/F}: \Bbb S_\infty\to \Bbb  C_\infty[G],$$
where $\Bbb S_\infty$ is a certain Goss space of mixed characteristic containing naturally a copy of $\Bbb Z_{\geq 0}$ and $\Bbb C_\infty$ is the completion of the algebraic closure of $k_\infty:=\Bbb F_q((t^{-1}))$ with respect to the unique extension of the valuation $v_\infty$ of $k$ of uniformizer $1/t.$
(See the Introduction in \cite{FGHP} for the precise definitions.) 

As a natural $G$--equivariant generalization of the Goss $\zeta$--function $\zeta_F^E$ associated to $E$ and studied by Taelman in \cite{Taelman}, the $L$--function $\Theta_{K/F}^E$ is defined as an infinite product of Euler factors 
$$\Theta_{K/F}^E(s):={\prod_{v}}'P_v^{\ast, G}(Nv^{-s})^{-1},$$
where $v$ runs over all the primes in the maximal spectrum ${\rm MSpec}(\mathcal O_F)$ of $\mathcal O_F$ which are tamely ramified in $K/F$ and of good reduction for $E$. For each such prime $v$, $P_v^{\ast, G}(X)$ is a polynomial with coefficients in $A[G]$, very closely related to the well understood characteristic polynomial
$$P_v(X):={\rm det}_{A_{v_0}}(X\cdot I_r-\sigma(v)\mid T_{v_0}(E))$$
of the action of a Frobenius morphism $\sigma(v)\in{\rm Gal}(\overline F/F)$ associated to $v$ on the $v_0$--adic Tate module $T_{v_0}(E)$ associated to the rank $r$ Drinfeld module $E$ at any prime $v_0\in{\rm MSpec}(A)\setminus\{v\cap A\}$, viewed as a free, rank $r$ module over the completion $A_{v_0}$ of $A$ at $v_0$. In fact, one has equalities (see \cite{FGHP}, Introduction):
$$P_v^{\ast, G}(1)=\frac{P_v(e_v\sigma_v)}{P_v(0)}, \qquad \Theta_{K/F}^E(0)={\prod_v}' \frac{P_v(0)}{P_v(e_v\sigma_v)},$$
where the infinite product converges in $\Bbb F_q((t^{-1}))[G]$, $\sigma_v\in G$ is a choice of Frobenius for $v$ in $G$, and $e_v$ is the idempotent in $A[G]$ associated to the trivial character of the inertia group of $v$ in $G$.
\medskip

The main result of \cite{FGHP} is the proof of a Tamagawa number formula of the type
$$\Theta_{K/F}^{E, \mathcal M}(0)=\frac{{\rm Vol}(E(K_\infty/\mathcal M))}{{\rm Vol}(K_\infty/\mathcal M)},$$
where $\Theta_{K/F}^{E, \mathcal M}(0)$ is a certain Euler product completion of 
$\Theta_{K/F}^E(0)$ at primes $v\in{\rm MSpec}(\mathcal O_F)$ which are either wildly ramified in $K/F$ or of bad reduction for $E$, constructed out of some additional arithmetic data $\mathcal M$ (called a taming module) for $K/F$. The numerator and denominator of the right side in the above formula are both volumes (with values in $\Bbb F_q((t^{-1}))[G]$ and defined precisely in \cite{FGHP}) of certain compact $A[G]$--modules of Arakelov type $E(K_\infty/\mathcal M)$ and $K_\infty/\mathcal M$.
\medskip

The above formula, proved in \cite{FGHP} for Drinfeld modules and extended in \cite{Green-Popescu} for higher dimensional abelian $t$--modules, generalizes to the $G$--equivariant setting Taelman's celebrated class--number formula \cite{Taelman} for the value $\zeta_F^E(0)$ of the Goss zeta--function associated to $E$. The interested reader should also consult \cite{Beaumont} for a slightly different proof in the abelian $t$--module case, using deformation theory, which only gives the desired result under certain restrictive conditions.
\medskip

However, the proof of the above formula in \cite{FGHP} and \cite{Green-Popescu} (and the same applies to the proofs given in \cite{Taelman} and \cite{Beaumont}) hinges upon an essential equality at the level of Euler factors, namely
$$\frac{P_v(e_v\sigma_v)}{P_v(0)}=\frac{|E(\mathcal O_K/v)|_G}{|\mathcal O_K/v|_G}\qquad \text{ in } \Bbb F_q((t^{-1}))[G],$$
for all good and tame primes $v$ as above, where the numerator and denominator of the right side  
are certain special generators of the Fitting ideals of the finite $A[G]$--modules of finite projective dimension $E(\mathcal O_K/v)$ and $\mathcal O_K/v$, respectively. 
In the Appendix of \cite{FGHP}, we prove the above equality only in the simplest case where $E:=\mathcal C$ is the Carlitz module, and Taelman does the same in \cite{Taelman}, in the case where $G$ is trivial. 
\medskip

In this paper, we develop techniques which allow us to prove the above equality for arbitrary Drinfeld modules (see Theorem \ref{main-theorem} and its proofs in the unramified case, given in \S\ref{unramified-section} and tamely ramified case, given in \S\ref{tame-section}.) In \S\ref{tmodule-section}, we also indicate how our techniques can be easily extended to prove the above formula for higher dimensional abelian $t$--modules, satisfying a certain purity condition. The case of arbitrary abelian $t$--modules will be treated in a separate, upcoming paper.
\medskip

{\bf Acknowledgement.} We would like to thank Urs Hartl for his expert and generous help with the material on local shtukas in \S\ref{motive-section}. His private e-mail exchanges with the first author \cite{private} were invaluable to us.

\section{The statement of the problem}\label{statement-section}
Let $p$ be a prime number and let $q$ be a power of $p$. In what follows, $k:=\mathbb{F}_q(t)$ denotes the rational function field in one variable over $\mathbb{F}_q.$ For any commutative $\mathbb{F}_q$-algebra $R,$ we denote by $\tau$ the $q$-power Frobenius endomorphism of $R$. We denote by $R\{\tau\}$ the twisted polynomial ring in $\tau$, with the property that
$$\tau\cdot x = x^q\cdot\tau \hspace{10pt} \forall \hspace{5pt} x \in R.$$

Let $F$ be a finite, separable extension of $\mathbb{F}_q(t)$ and let $K$ be a finite abelian extension of $F$ with Galois group $G.$ We also assume that the field of constants in $K$ is $\mathbb{F}_q,$ i.e.
$$K \cap \overline{\mathbb{F}}_q = \mathbb{F}_q.$$
Let us denote $\mathbb{F}_q[t]$ by $A.$ Note that if $v$ denotes an arbitrary normalized valuation on $\Bbb F_q(t)$ and $\infty$ denotes the normalized valuation of uniformizer $1/t$, then
$$A = \{a \in \mathbb{F}_q(t) \,\big|\, v(a) \geq 0, \text{ for all } \hspace{5pt} v \neq \infty\}.$$
Let $\mathcal{O}_F$ and $\mathcal{O}_K$ denote the integral closures of $A$ in $F$ and $K$, respectively. In what follows, we abuse notation and use the same letter for normalized valuations and the associated maximal ideals of elements of strictly positive valuation. 
\medskip 

Next, we consider a Drinfeld module $E$ of rank $r \in \mathbb{N}$ defined on $A$ with values in $\mathcal O_F\{\tau\}.$ More precisely, $E$ is given by an $\mathbb{F}_q$-algebra morphism
$$\phi_E: A \to \mathcal{O}_F\{\tau\}, \hspace{20pt} t \mapsto t\cdot\tau^0 + e_1\tau+\dots + e_r\cdot\tau^r,$$
where $a_i \in \mathcal{O}_F$, for all $i$  and $e_r \neq 0.$ This gives rise to  a functor
$$E: (\mathcal{O}_F\{\tau\}[G]-\text{modules}) \to (A[G]-\text{modules}).$$
In other words, for any $\mathcal{O}_F\{\tau\}[G]$-module $M,$ we denote by $E(M)$ the $A[G]-$module whose underlying $\Bbb F_q[G]$--module is $M$ and the $A$-action is given by
$$t \star m = \phi_E(t)\cdot m = t\cdot m + e_1\tau\cdot m + \dots  + e_r\tau^r\cdot m$$
\medskip

Let $v_0 \in {\rm MSpec}(A)$ and let $A_{v_0}$ and $k_{v_0}$ denote the completions of $A$ and $k$ with respect to the valuation $v_0$. For all $n \in \mathbb{N}$, we denote by $E[v_0^n]$ the $A_{v_0}$-module of $v_0^n$-torsion points of $E,$ i.e.
$$E[v_0^n] = \{x \in E(\overline{F})| f \star x = 0, \text{ for all }f \in v_0^n\}.$$
The $v_0$-adic Tate module of $E$ is defined as
$$T_{v_0}(E) = \text{Hom}_{A_{v_0}}(k_{v_0}/A_{v_0}, E[v_0^\infty]).$$
Since $A$ is a PID, we also have
$$T_{v_0}(E) = \varprojlim E[v_0^n],$$
where the transition maps in the projective limit are given by multiplication with a generator of $v_0$, while  $E[v_0^\infty] = \bigcup_{n \geq 1} E[v_0^n].$ 
Recall that $E[v_0^n]$ and $T_{v_0}(E)$ are free modules of rank $r$ over $A/v_0^n$ and $A_{v_0}$, respectively, and are endowed with obvious $A_{v_0}$-linear, continuous $G_F$-actions, where $G_F = Gal(\overline{F}/F).$\\

Let $v \in$ MSpec($\mathcal{O}_F$), such that $v \nmid v_0.$ Fix a choice of decomposition group $G(v) \subset G_F,$ and a Frobenius morphism $\sigma(v) \in G(v).$ Then, it is known (see \cite{FGHP} and the references therein) that if $E$ has good reduction at $v$ (i.e. $v\nmid e_r$), the $G_F$-representation $T_{v_0}(E)$ is unramified at $v$ and the polynomial
$$P_v(X) = {\rm det}_{A_{v_0}}(X\cdot I_r - {\sigma(v)}|T_{v_0}(E))$$
is independent of $v_0$ and actually lies in $A[X].$ Above, $I_r$ denotes the $r\times r$ identity matrix.

\begin{Def}\label{monic-gen-def}
    Let $M$ be an $A[G]$-module which is free of rank $m$ as an $\mathbb{F}_q[G]$-module. Then it is known (see \cite{FGHP} Proposition A.4.1) that the Fitting ideal ${\rm Fitt}_{A[G]}^{0}(M)$ is principal and has a unique $t$--monic generator $f_M(t) \in A[G] = \mathbb{F}_q[G][t]$ of degree $m.$ We denote this generator by  $|M|_G,$ i.e. 
    $$|M|_G = f_M(t) \in \mathbb{F}_q[G][t].$$
\end{Def}

The following is Proposition A5.1. from the Appendix in \cite{FGHP}:
\begin{proposition}\label{freeness-prop}
    Assume that $v$ is tamely ramified in $K/F$ and let $E$ be any Drinfeld module as above. Let $w_0$ denote the prime in $A$ sitting below $v$ and let $f(v/w_0) = [\mathcal{O}_F/v:A/w_0].$ Then the following hold:
    \begin{enumerate}
        \item The $\mathbb{F}_q[G]$-modules $\mathcal{O}_K/v$ and $E(\mathcal{O}_K/v)$ are free of rank $n_v = [\mathcal{O}_F/v:\mathbb{F}_q]$ and therefore $|\mathcal O_K/v|_G$ and $|E(\mathcal O_K/v)|_G$ are monic polynomials of $t$--degree $n_v$.

        \item We have an equality
        $$|\mathcal{O}_K/v|_G = Nv$$
        where $Nv$ denotes the unique monic generator of $w_0^{f(v/w_0)}$ and $f(v/w_0):=[\mathcal O_F/v:A/w_0]$.
    \end{enumerate}
\end{proposition}

Let $I_v \subset G_v \subset G$ denote the inertia and decomposition groups of $v$ in $G$, respectively. Let $\sigma_v$ be the image of ${\sigma(v)}$ via the Galois restriction map $G(v)\twoheadrightarrow G_v$. Our main goal in this paper is the proof of the following.
\begin{theorem}\label{main-theorem}
Assume that $v \in{\rm MSpec}(\mathcal{O}_F) $ is tamely ramified in $K/F$ and that $E$ has good reduction at $v$. Then, we have an equality in $\Bbb F_q[G][[1/t]]$
$$\frac{P_v(\sigma_ve_v)}{P_v(0)} = \frac{|E(\mathcal{O}_K/v)|_G}{|\mathcal{O}_K/v|_G},$$
where $e_v = \frac{1}{|I_v|}\sum_{\sigma \in I_v}\sigma$ is the idempotent of the trivial character of $I_v$ in $A[G].$
\end{theorem}

A proof of the above statement in the case where $E$ is the Carlitz module $C$, defined by $\phi_C(t)=t+\tau$, was given in the Appendix of \cite{FGHP}. Below, we develop techniques which settle the above theorem for a general Drinfeld module $E$.
Proposition 1.2(2) gives us a good understanding of $|\mathcal{O}_K/v|_G$. Therefore a major portion of our work is directed towards understanding the relation between $|E(\mathcal{O}_K/v)|_G$ and $P_v(\sigma_ve_v).$ \\

\section{The reduction $\overline E$ of $E$ modulo $v$}\label{reduction-section}
    In this section, we fix a prime $v\in\text{MSpec}(\mathcal O_F)$ such that $E$ has good reduction at $v.$ We are not assuming that $v$ is necessarily tamely ramified in $K/F$. Let us denote by $w_0$ the prime in $A$ that lies below $v.$ After reduction of $E \mod v$ (i.e. reduction of the coefficients of $E$ modulo $v$), we obtain the rank $r$ Drinfeld module $\overline{E},$ defined over $\mathcal O_F/v$, given by the $\mathbb{F}_q$-algebra morphism
    $$\phi_{\overline{E}}: A \to\mathcal O_F\{\tau\}\twoheadrightarrow \mathcal{O}_F/v\{\tau\}$$
    where $\phi_{\overline{E}}(t) = i(t)\cdot \tau^0 + ... + i(e_r)\cdot \tau^r$ with $i: A \xhookrightarrow{} \mathcal{O}_F \twoheadrightarrow \mathcal{O}_F/v$ being the obvious map. The Drinfeld module $\overline E$ has rank $r$, characteristic $w_0:=\ker(i)$, and height $h$. By definition, $h$ is the unique integer $0<h\leq r$ which, for any generator $\pi_{w_0}$ of $w_0$, satisfies the equality
    $$\phi_{\overline E}(\pi_{w_0})=b_{d_0h}\tau^{d_0h}+b_{d_0h+1}\tau^{d_0h+1}+\dots+b_{d_0r}\tau^{rd_0},$$
    where $d_0:=[A/w_0:\Bbb F_q]$, $b_i\in\mathcal O_F/v$, and $b_{d_0h}\ne 0$. (See \cite{Goss}, Section 4.5 for the existence of $h$.)
    \medskip
    
    Recall that, by the notation introduced above, we have a field isomorphism $\mathcal{O}_F/v \simeq \mathbb{F}_{q^{n_v}}.$ Consequently, the Frobenius morphism ${\rm Frob}_{q^{n_v}}:=\tau^{n_v}$ is an endomorphism of $\overline E$ and it acts naturally and linearly on all the Tate modules associated to $\overline E$. \\ 

    Next, we fix $v_0\in\text{MSpec}(A)$, $v_0\ne w_0$ and consider the characteristic polynomial of the action of the $q^{n_v}$-power Frobenius morphism, viewed as an endomorphism of $\overline E$,  on the free $A_{v_0}$--module $T_{v_0}(\overline{E})$ of rank $r$:
    $$f_{\overline{E}}(X) = {\rm det}_{A_{v_0}}(X\cdot I_r - {\rm Frob}_{q^{n_v}} \mid T_{v_0}(\overline{E})).$$
    Then, $f_{\overline{E}}(X)$ is independent of $v_0$ and lies in $A[X]$. (See \S4.12 in \cite{Goss}.) By Theorem 4.12.15 in \cite{Goss} and the discussion preceding that, we have the following.
    \begin{proposition}\label{eigenvalues-prop}
        Any root $\alpha$ of $f_{\overline E}$ satisfies the following properties:
        \begin{enumerate}
            \item $w(\alpha) = 0$ for all finite places $w$ of $\mathbb{F}_q(t)(\alpha),$ except for exactly one place above $w_0.$

            \item There is only one place of $\mathbb{F}_q(t)(\alpha)$ lying above $\infty.$

            \item $|\alpha|_\infty = q^{\frac{n_v}{r}}$ where $|.|_\infty$ denotes the unique extension to $\mathbb{F}_q(t)(\alpha)$ of the normalized absolute value of $\mathbb{F}_q(t)$ corresponding to $\infty.$

            \item $[\mathbb{F}_q(t)(\alpha):\mathbb{F}_q(t)]$ divides $r.$
        \end{enumerate}
    \end{proposition}
\medskip
    Further, one can consider the characteristic polynomial of ${\rm Frob}_{q^{n_v}}$ (viewed as endomorphism of $\overline E$) acting on the free $A_{w_0}$--module $T_{w_0}(\overline E)$ of rank $(r-h)$ (see \cite{Goss}, Section 4.5 for the calculation of the rank)
    $$g_{\overline{E}}(X) = {\rm det}_{A_{w_0}}(X\cdot I_{r-h} - {\rm Frob}_{q^{n_v}} \mid T_{w_0}(\overline{E})).$$
    This is a monic polynomial in $A_{w_0}[X]$ of degree $(r-h).$  In \S? below we will prove the following.

    \begin{proposition}\label{g-divides-f-prop} The polynomial $g_{\overline E}(X)$ divides the polynomial $f_{\overline E}(X)$ in $A_{w_0}[X]$.
    \end{proposition}

\begin{proof} See \S\ref{motive-section} below for the proof of a stronger statement.
\end{proof}

\medskip
    
    Let $\mathcal O_v$ and $F_v$ be the completions at $v$ of $\mathcal O_F$ and $F$, respectively. Our choice of decomposition group $G(v)$ corresponds to choosing an embedding  $\overline F\to \overline{F_v}$ at the level of separable closures of $F$ and $F_v$, such that Galois restriction induces a group isomorphism $G(\overline{F_v}/F_v)\simeq G(v)$. Since $E$ has good reduction at $v$ and the Galois representations $E[v_0^n]$ are unramified at $v$, it is not difficult to see that we have 
    $$E[v_0^n]\subseteq \mathcal O_v^{nr},\text{ for all }n\geq 1,$$ 
    where $\mathcal O_v^{nr}$ is the integral closure of $\mathcal O_v$ in the maximal unramified extension $F_v^{unr}$ of $F_v$ in $\overline{F_v}$. Moreover, the reduction $\mod v$ map induces isomorphisms of $A_{v_0}[[\overline {G(v)}]]$--modules
    \begin{equation}\label{Tate-module-iso} E[v_0^n]\simeq \overline E[v_0^n], \qquad T_{v_0}(E)\simeq T_{v_0}(\overline E),\end{equation}
    where 
    $$\overline{G(v)}:=G(v)/I(v)\simeq G(\overline{\Bbb F_{q^{n_v}}}/\Bbb F_{q^{n_v}}).$$
The group isomorphism above sends $\overline{{\sigma(v)}}$ (the image of our choice of Frobenius $\sigma(v)$ in $\overline {G(v)}$) to ${\rm Frob}_{q^{n_v}}$. Consequently, we have an equality of characteristic polynomials in $A[X]$:
   \begin{equation}\label{char-poly-equality} f_{\overline{E}}(X) = P_v(X).\end{equation}
Consequently, Proposition \ref{eigenvalues-prop} gives us information on the roots of the characteristic polynomial $P_v(X)$. The following corollary regarding the coefficients of $P_v(X)$ will be particularly useful in what follows.

\begin{Cor}\label{coeff-charpoly-cor} Let $P_v(X)=a_0+a_1 X+\cdots + a_{r-1}X^{r-1}+X^r$, with $a_0,\dots, a_{r-1}\in A$. Then, we have
\begin{enumerate}
\item ${\rm deg}_t(a_0)=n_v$ and $0<{\rm deg}_t(a_i)<n_v$, for all $i>0$.
\item $P_v(X)\in\Bbb F_q[X][t]$ is a polynomial of degree $n_v$ in $t$ with the same leading coefficient as $a_0$.
\item $a_0=\rho\cdot Nv$, for some $\rho\in\Bbb F_q^\times$, where $Nv$ is the unique monic generator of $w_0^{f(v/w_0)}.$
\end{enumerate}
Above, ${\rm deg}_t(\ast)$ denotes the degree in $t$ of a polynomial in $A=\Bbb F_q[t].$
\end{Cor}
\begin{proof}Let $\alpha_1, ..., \alpha_r \in \overline{A}$ denote the roots of $P_v(X)$ in the integral closure of $A$. Then
    
    $$\begin{aligned}
    P_v(X) = \prod_{i = 1}^{r}(X - \alpha_i) &= (-1)^r\prod_{i = 1}^{r}\alpha_i + (-1)^{r - 1}\bigg(\sum_{j = 1}^{r}\prod_{i \neq j}\alpha_i\bigg)X + ... + X^r\\
    &=a_0+a_1\cdot X+\dots+a_{r-1}\cdot X^{r-1}+X^r.\end{aligned}$$
Let $|\cdot|_\infty$ denote an extension to $\mathbb{F}_q(t)(\alpha_1, ..., \alpha_r)$ of the normalized absolute value of $\mathbb{F}_q(t)$ corresponding to $\infty$ (also denoted by $|\cdot|_\infty$ below.) By Proposition 1.4, we have $|\alpha_i|_\infty = q^{\frac{n_v}{r}}$, for all $i$. Therefore, we have
    $$|a_0|_\infty=\bigg\lvert\prod_{i = 1}^{r}\alpha_i\bigg\rvert_\infty = q^{n_v}, \qquad {\rm deg}_t(a_0)= \log_q(|a_0|_\infty)=n_v.$$
    Furthermore, since $|\cdot|_\infty$ is non-archimedean, we have 
    $$|a_i|_\infty\leq q^{n_v\cdot\frac{r-i}{r}},\qquad  {\rm deg}_t(a_i)= \log_q(|a_i|_\infty)\leq n_v\cdot\frac{r-i}{r}<n_v, \quad\text{ for all }i\geq 1. $$
    This concludes the proof of part (1). 
    \medskip
    
    Part (2) is a direct consequence of part (1). Further, since $a_0=(-1)^r\prod_{i=1}^r\alpha_i$, part (3) is a direct consequence of Proposition \ref{eigenvalues-prop}(1)--(3).
\end{proof}

\bigskip

\section{Fitting ideals of Tate modules and consequences}\label{Fitting-Tate-section}   

    We begin by stating a general commutative algebra result regarding Fitting ideals of modules over certain rings of equivariant Iwasawa algebra type. For a proof of this result, see Proposition 4.1 in \cite{GP}.
    \begin{proposition}[Greither--Popescu]
        Let $R$ be a semi-local, compact topological ring, and let $\Gamma$ be a pro-cyclic group, topologically generated by $\gamma.$ Suppose that $M$ is an $R[[\Gamma]]$--module which is free of rank $n$ as an $R$--module. Let $M_\gamma \in M_n(R)$ denote the matrix of the action of $\gamma$ on some $R$-basis of $M.$ Then, we have an equality of $R[[\Gamma]]$--ideals
        $${\rm Fitt}_{R[[\Gamma]]}(M) = \bigg({\rm det}_R (X\cdot I_n - M_\gamma)\bigg\lvert_{X = \gamma}\bigg)$$
    \end{proposition}
\begin{proof} See the proof of Proposition 4.1 in \cite{GP}.
\end{proof}
    
An immediate consequence of the above proposition is the following.
\begin{Cor}\label{Fitting-Tate-corollary}
For all $v_0\in{\rm MSpec}(A)$, we have the following equalities of $A_{v_0}[[\overline{G(v)}]]$--ideals.
\begin{enumerate}
\item If $v_0\ne w_0$, then 
$${\rm Fitt}_{A_{v_0}[[\overline{G(v)}]]}(T_{v_0}(\overline E))=
{\rm Fitt}_{A_{v_0}[[\overline{G(v)}]]}(T_{v_0}(E)) = (f_{\overline E}(\overline{\sigma(v)})) = (P_v(\overline{\sigma(v)})).$$
\item If $v_0=w_0$, then 
$${\rm Fitt}_{A_{w_0}[[\overline{G(v)}]]}(T_{w_0}(\overline E))=(g_{\overline E}(\overline{\sigma(v)}))$$
\end{enumerate}
\end{Cor}

\begin{proof}
Apply the proposition above to $R:=A_{v_0}$, $\Gamma:=\overline{G(v)}$, $\gamma:=\overline{\sigma(v)}$ and the module 
$$M:=T_{v_0}(E)\simeq T_{v_0}(\overline E),$$
which is $A_{v_0}$--free of rank $r$, if $v_0\ne w_0$ and, respectively, the module
$$M:=T_{w_0}(\overline E),$$
which is $A_{w_0}$--free of rank $(r-h)$, if $v_0=w_0.$
\end{proof}

Next, we fix a prime $w\in{\rm MSpec}(\mathcal O_K)$ lying above $v$. We let $G(w)=G(v)\cap G_K$, $I(w)=I(v)\cap G_K$ and $\sigma(w):=\sigma(v)^{f}$, where $f: =f(w/v)=[\mathcal O_K/w:\mathcal O_F/v]$. Then, $\sigma(w)\in G(w)$ is a choice of Frobenius for $w$ and its image $\overline{\sigma(w)}\in \overline {G(w)}:=G(w)/I(w)$ corresponds via the group isomorphism $\overline{G(w)}\simeq G_{\mathcal O_K/w}$ to ${\rm Frob}_{q^{fn_v}}$, viewed as an endomorphism of $\overline E$.

\begin{Lemma}\label{coinvariants-Tate-lemma} Let $v_0\in{\rm MSpec}(A)$. Then, we have the following canonical isomorphisms of $A_{v_0}[G_v]$--modules:
\begin{enumerate}
\medskip
\item     
    $T_{v_0}(\overline{E})/(1 - \overline{\sigma(w)})T_{v_0}(\overline{E}) \simeq \overline{E}(\mathcal{O}_K/w)_{v_0}={E}(\mathcal{O}_K/w)_{v_0}.$
\medskip
\item 
    $T_{v_0}(E)/(1 - \sigma(w))T_{v_0}(E) \simeq E(\mathcal{O}_K/w)_{v_0}$, assuming that $v_0\ne w_0$.
\end{enumerate}
\medskip
    Above, we let $M_{v_0}:=M\otimes_AA_{v_0}$ for any $A$--module $M$.
    \end{Lemma}
    \begin{proof} By the definition of $\overline E$, we have 
    $$E(\mathcal O_K/w)=\overline E(\mathcal O_K/w).$$ 
    Consequently, part (2) follows from part (1) via the second isomorphism in \eqref{Tate-module-iso} above. 
    
    In order to prove part (1), apply the functor $\ast\to {\rm Hom}_{A_{v_0}}(\ast, \overline E[v_0^\infty])$ to the exact sequence of $A_{v_0}$--modules
    $$0\longrightarrow A_{v_0}\longrightarrow k_{v_0}\longrightarrow k_{v_0}/A_{v_0}\longrightarrow 0.$$
    Since the $A_{v_0}$--module $\overline E[v_0^\infty]:=\cup_n\overline E[v_0^n]$ is divisible and therefore injective (as $A_{v_0}$ is a PID), the above functor is exact. Therefore, we obtain the following exact sequence of $A_{v_0}[\overline{G(v)}]$--modules:
    \begin{equation}\label{exact-sequence-Tate}0\longrightarrow T_{v_0}(\overline E)\longrightarrow {\rm Hom}_{A_{v_0}}(k_{v_0}, \overline E[v_0^\infty])\longrightarrow \overline E[v_0^\infty]\longrightarrow 0.\end{equation}
    Now, it is easy to see that one has an isomorphism of $k_{v_0}[\overline{G(v)}]$--modules
    $$k_{v_0}\otimes_{A_{v_0}}T_{v_0}(\overline E)\simeq {\rm Hom}_{A_{v_0}}(k_{v_0}, \overline E[v_0^\infty]), \quad \xi\otimes\phi \to (x\to \phi(\widehat{\xi\cdot x})),$$
    for all $\xi, x\in k_{v_0}$ and all $\phi\in T_{v_0}(\overline E)={\rm Hom}_{A_{v_0}}(k_{v_0}/A_{v_0}, \overline E[v_0^\infty])$, where $\widehat{x\cdot\xi}$ is the class of $x\cdot \xi$ in $k_{v_0}/A_{v_0}$.

    Now, Proposition \ref{eigenvalues-prop}(3) and Proposition \ref{g-divides-f-prop} show that the eigenvalues of $\overline{\sigma(w)}=\overline{\sigma(v)}^f=({\rm Frob}_{q^{n_v}})^f$ acting on the $k_{v_0}$--vector space $k_{v_0}\otimes_{A_{v_0}}T_{v_0}(\overline E)$ are all different from $1$. Consequently, $(\overline{\sigma(w)}-1)$ is an automorphism of this $k_{v_0}$--vector space. Consequently, when one takes the $\overline{\sigma(w)}$--invariants and coinvariants in the exact sequence \eqref{exact-sequence-Tate} above, one obtains an isomorphism of $A_{v_0}[G_v]$--modules
    $$T_{v_0}(\overline{E})/(1 - \overline{\sigma(w)})T_{v_0}(\overline{E}) \simeq \overline E[v_0^\infty]^{\overline{\sigma(w)}=1}=\overline{E}(\mathcal{O}_K/w)_{v_0},$$
    which concludes the proof of the Lemma.
    
    \end{proof}
    \begin{Cor}\label{Fitting-residuefield-cor} For all $v_0\in{\rm MSpec}(A)$, the following equalities of $A_{v_0}[\overline{G_v}]$--ideals hold:
    $${\rm Fitt}_{A_{v_0}[\overline{G_v}]}E(\mathcal O_K/w)_{v_0}=\begin{cases} \big(P_v(\overline{\sigma_v}\big)=\big(f_{\overline E}(\overline{\sigma_v})\big), & \text{ if }v_0\ne w_0\\
    \big(g_{\overline E}(\overline{\sigma_v})\big), & \text{ if } v_0=w_0.
    \end{cases}
    $$
    Further, if $v_0\ne w_0$, then we have an equality of $A_{(v_0)}[\overline{G_v}]$--ideals
    $${\rm Fitt}_{A_{(v_0)}[\overline{G_v}]}E(\mathcal O_K/w)_{v_0}=\big(P_v(\overline{\sigma_v}\big)=\big(f_{\overline E}(\overline{\sigma_v})\big).$$
    Here, $\overline G_v:=G_v/I_v$, $\overline\sigma_v$ is the image of $\sigma_v$ in $\overline{G_v}$, and $A_{(v_0)}$ is the localization of $A$ at $v_0$.
    \end{Cor}
    \begin{proof} First, note that the isomorphism of $A_{v_0}[G_v]$--modules in Lemma \ref{coinvariants-Tate-lemma}(1) can be rewritten as an isomorphism of $A_{v_0}[\overline{G_v}]$--modules
    $$T_{v_0}(\overline E)\otimes_{A_{v_0}[[\overline{G_v}]]} A_{v_0}[\overline{G(v)}]\simeq E(\mathcal O_K/w)_{v_0},$$ where the ring morphism $\pi: A_{v_0}[[\overline {G(v)}]]\twoheadrightarrow A_{v_0}[\overline{G_v}]$ is the $A_{v_0}$--linear map given by Galois restriction, which maps $\overline{\sigma(v)}\to\overline{\sigma_v}.$ The isomorphism above permits us to apply the well known base--change property of Fitting ideals which, combined with Corollary \ref{Fitting-Tate-corollary} above, leads to the equalities of $A_{v_0}[\overline{G_v}]$--ideals
    $${\rm Fitt}_{A_{v_0}[\overline{G_v}]}E(\mathcal O_K/w)_{v_0}=\pi\big({\rm Fitt}_{A_{v_0}[[\overline{G(v)}]]}(T_{v_0}(E))\big)=\begin{cases} \big(P_v(\overline{\sigma_v}\big)=\big(f_{\overline E}(\overline{\sigma_v})\big), & \text{ if }v_0\ne w_0\\
    \big(g_{\overline E}(\overline{\sigma_v})\big), & \text{ if } v_0=w_0.
    \end{cases}.$$
    Next, assume that $v_0\ne w_0$ and observe that since $E(\mathcal O_K/w)$ is finite (and therefore $A$--torsion), we have isomorphisms 
    $$E(\mathcal O_K/w)\otimes_AA_{(v_0)}\simeq \big(E(\mathcal O_K/w)\otimes_AA_{(v_0)}\big)\otimes_{A_{(v_0)}} A_{v_0}\simeq E(\mathcal O_K/w)_{v_0}.$$
    Consequently, base--change for Fitting ideals applied to the ring extension $A_{(v_0)}[G]\subseteq A_{v_0}[G]$ and the last equality of ideals displayed above gives
    $$P_v(\overline{\sigma_v})A_{v_0}[G]={\rm Fitt}_{A_{(v_0)}[\overline{G_v}]}\big(E(\mathcal O_K/w)_{v_0}\big)A_{v_0}[\overline{G_v}].$$
    However, since the ring extension $A_{(v_0)}[G]\subseteq A_{v_0}[G]$ is faithfully flat (because $A_{(v_0)}\subseteq A_{v_0}$ is), we have 
    $$\begin{aligned}{\rm Fitt}_{A_{(v_0)}[\overline{G_v}]}\big(E(\mathcal O_K/w)_{v_0}\big)&={\rm Fitt}_{A_{(v_0)}[\overline{G_v}]}\big(E(\mathcal O_K/w)_{v_0}\big)A_{v_0}[\overline{G_v}]\cap A_{(v_0)}[\overline{G_v}]\\ 
    &=P_v(\overline{\sigma_v})A_{v_0}[\overline{G_v}]\cap A_{(v_0)}[\overline{G_v}]=P_v(\overline{\sigma_v})A_{(v_0)}[\overline{G_v}].
    \end{aligned}$$
    Above, we used the fact that if $R\subseteq R'$ is a faithfully flat extension of commutative rings and $I$ is an ideal in $R$, then $IR'\cap R=I$. (See \cite{Matsumura}, Chapter 2, Section 4, 4.C(ii).) 
    \end{proof}

In the next two sections, we provide a proof of Theorem \ref{main-theorem}. For technical reasons which will become apparent shortly, we treat first the unramified case.
    
\section{The unramified case}\label{unramified-section}
    We keep the notations and assumptions of the previous section. In addition, we assume  that the prime $v$ is unramified in $K/F.$ Consequently, we have $\overline{G_v}= G_v$ and $\overline{\sigma_v}=\sigma_v$ throughout.  
    
    \begin{Lemma} \label{fitt-local-lemma}Under the current assumptions, we have equalities of ideals.
      \begin{enumerate}
          \item 
      ${\rm Fitt}_{A_{(v_0)}[G]}(E(\mathcal{O}_K/v)_{v_0}) = (P_v(\sigma_v))=(f_{\overline E}(\sigma_v)),$ 
      for all $v_0\in{\rm MSpec}(A)$, with $v_0\ne w_0$.
    \item ${\rm Fitt}_{A_{w_0}[G]}(E(\mathcal{O}_K/v)_{w_0}) = (g_{\overline E}(\sigma_v)).$
    \end{enumerate}
    \end{Lemma}
    \begin{proof} In this case, we have an isomorphism of $A_{(v_0)}[G]$--modules
    $$E(\mathcal O_K/v)_{v_0}\simeq E(\mathcal O_K/w)_{v_0}\otimes_{A_{(v_0)}[G_v]}A_{(v_0)}[G],$$
for all $v_0\in{\rm MSpec}(A)$. (See the Appendix of \cite{FGHP}.) Therefore, the equality in the Lemma follows from Corollary \ref{Fitting-residuefield-cor} and the base--change property of Fitting ideals.
    \end{proof}

    We are ready to prove the following refinement of Theorem \ref{main-theorem} in the unramified case. 

    \begin{proposition}\label{prop-unr} Assume that $v$ is unramified in $K/F$. Let $\rho\in\Bbb F_q^\times$ as defined in Corollary \ref{coeff-charpoly-cor}(3). Then, $\rho^{-1}P_v(0)$ and $\rho^{-1}P_v(\sigma_v)$ are monic polynomials of $t$--degree $n_v$ in $A=\Bbb F_q[t]$ and $A[G]=\Bbb F_q[G][t]$, respectively, and the following hold.
    \begin{enumerate}
\item $\rho^{-1}P_v(0)=|\mathcal O_K/v|_G.$
\item $\rho^{-1}P_v(\sigma_v)=|E(\mathcal O_K/v)|_G.$
\item $P_v(\sigma_v)/P_v(0)= |E(\mathcal O_K/v)|_G/|\mathcal O_K/v|_G$ in $\Bbb F_q[G][[t]]$, i.e. Theorem \ref{main-theorem} holds.
    \end{enumerate}
        
    \end{proposition}

    \begin{proof}
    According to Corollary \ref{coeff-charpoly-cor} (see (1) and (3) in loc.cit.), $P_v(0)\in \Bbb F_q[t]$ and $P_v(X)\in \Bbb F_q[X][t]=A[X]$,
    viewed as polynomials in $t$, have degrees equal to $n_v$ and leading coefficients equal to $\rho$. Therefore, $\rho^{-1}P_v(0)$ and $\rho^{-1}P_v(\sigma_v)$ are indeed monic polynomials of common $t$--degree $n_v$. Further, Corollary \ref{coeff-charpoly-cor}(3) shows that $\rho^{-1}P_v(0)=Nv$ which, if combined to Proposition \ref{freeness-prop}(2), proves part (1) of the statement above.
    \medskip

Next, we focus on the proof of equality (2) in the Proposition above. In order to simplify the notation, let 
$$f:=\rho^{-1}P_v(\sigma_v), \qquad g:=\big|E(\mathcal O_K/v)\big|_G.$$
Then, $f$ and $g$ are both monic polynomials in $t$, of degrees equal to $n_v$. (See Proposition \ref{freeness-prop}(2) and Corollary \ref{coeff-charpoly-cor}.) Further,
Lemma \ref{fitt-local-lemma}(1) and the definition of $g$ imply that they satisfy the following equalities
$$fA_{(v_0)}[G]=gA_{(v_0)}[G]={\rm Fitt}_{A_{(v_0)}[G]}E(\mathcal O_K/v)_{v_0}, \quad \text{ for all }v_0\in\text{MSpec}(A)\setminus\{w_0\}.$$
Now, it is easy to check that the total ring of fractions of $A[G]$ is $k[G]$. Since $g$ and $f$ are monic, they are not zero--divisors in $A[G]$ and $k[G]$. Therefore, the equalities above imply that 
$$\frac{f}{g}\in k[G]\cap\big(\bigcap_{v_0\ne w_0}A_{(v_0)}[G]^\times\big)\subseteq \big(\bigcap_{v_0\ne w_0}A_{(v_0)}\big)[G].$$
This implies that there exists $\xi\in A[G]$ and $m\in\Bbb Z_{\geq 0}$, such that 
$$\frac{f}{g}=\frac{\xi}{\pi_{w_0}^m},\qquad \text{ with $\pi_{w_0}\nmid\xi$ in $A[G]$ if }m\geq 1,$$
where $\pi_{w_0}\in A$ is the unique monic generator of the maximal ideal $w_0$.

We claim that $m=0$. In order to prove that, let us note that Lemma \ref{fitt-local-lemma} implies that 
$$g_{\overline E}(\sigma_v)A_{w_0}[G]=gA_{w_0}[G].$$ 
On the other hand, Proposition \ref{g-divides-f-prop} implies that 
$$g_{\overline E}(\sigma_v)\mid f_{\overline E}(\sigma_v) \text{ in } A_{w_0}[G].$$
Consequently, since $f=\rho^{-1}P_v(\sigma_v)=\rho^{-1}f_{\overline E}(\sigma_v)$, we have
$$g\mid f \text{ in }A_{w_0}[G].$$
Therefore, $\xi/\pi_{w_0}^m\in A_{w_0}[G]$, which implies that 
$$\pi_{w_0}^m\mid \xi\text{ in }A_{w_0}[G].$$
However, since $\pi_{w_0}, \xi\in A[G]$, this shows that the divisibility above happens in $A[G]$. Therefore, we have $m=0$, as claimed and, consequently
$$\frac{f}{g}\in A[G]=\Bbb F_q[G][t].$$
However, since $f$ and $g$ are monic polynomials in $t$ of common degree $n_v$, we must have
$$f=g,$$
which concludes the proof of part (2) of the Proposition.
\medskip
Part (3) is a consequence of parts (1) and (2).

\end{proof}

 \begin{rem}\label{remark-unr}
Note that since Corollary \ref{coeff-charpoly-cor} is valid in general, regardless of the ramification status of $v$ in $K/F$, the polynomials $\rho^{-1}P_v(X)$, $\rho^{-1}P_v(e_v\sigma_v)$, and $\rho^{-1}P_v(0)$ are monic of degree $n_v$ in $t$ in general, even in the tamely ramified case. Also, the proof given to part (1) of the Proposition above is valid in general, so even in the tamely ramified case. 
\end{rem}

\section{The tamely ramified case}\label{tame-section}

Now, suppose that $v$ is tamely ramified in $K/F$, i.e. $p\nmid |I_v|$. As above, we denote by 
$$e_v:=\frac{1}{I_v}\sum_{\sigma\in I_v}\sigma$$
the idempotent in $A[G]$ associated to the trivial character of the inertia group 
$I_v$ of $v$ in $G$. The goal of this section is to prove the following analogue of Proposition \ref{prop-unr} in this case.
\begin{proposition}\label{prop-tame} Assume that $v$ is tamely ramified in $K/F$. Let $\rho\in\Bbb F_q^\times$ as defined in Corollary \ref{coeff-charpoly-cor}(3). Then, $\rho^{-1}P_v(0)$ and $\rho^{-1}P_v(e_v\sigma_v)$ are monic polynomials of $t$--degree $n_v$ in $A=\Bbb F_q[t]$ and $A[G]=\Bbb F_q[G][t]$, respectively, and the following hold.
    \begin{enumerate}
\item $\rho^{-1}P_v(0)=|\mathcal O_K/v|_G.$
\item $\rho^{-1}P_v(e_v\sigma_v)=|E(\mathcal O_K/v)|_G.$
\item $P_v(e_v\sigma_v)/P_v(0)= |E(\mathcal O_K/v)|_G/|\mathcal O_K/v|_G$ in $\Bbb F_q[G][[t]]$, i.e. Theorem \ref{main-theorem} holds.
    \end{enumerate}        
    \end{proposition}
\begin{proof}
Observe that part (3) is a consequence of parts (1) and (2). Since Remark \ref{remark-unr} settles everything but part (2) in the more general, tamely ramified case, we will focus below on proving part (2). 
\medskip

Note that, if $e:=|I_v|$ and $w$ is a fixed prime in $K$ above $v$, as before we have isomorphisms of $A[G]$--modules 
$$\mathcal O_{K}/v\simeq \mathcal O_K/w^e\otimes_{A[G_v]}A[G], \qquad E(\mathcal O_{K}/v)\simeq E(\mathcal O_K/w^e)\otimes_{A[G_v]}A[G].$$
Consequently, base--change for Fitting ideals and the isomorphism $\mathcal O_K/w^e\simeq \mathcal O_w/w^e$ show that it suffices for us to work locally and prove that we have an equality of $A[G_v]$--ideals 
\begin{equation}\label{Fitt-local}{\rm Fitt}_{A[G_v]}E(\mathcal O_w/w^e)=(P_v(e_v\sigma_v)),\end{equation}
where $\mathcal O_w$ is the ring of integers in the $w$--adic completion $K_w$ of $K$ and $G_v$ is now viewed as $G(K_w/F_v)$, where $F_v$, $\mathcal O_v$ etc. have the obvious meanings. Also, let us note that the faithful flatness of the ring extension $A[G_v]\subseteq A[G]$ combined with part (1) of the Proposition give us
$${\rm Fitt}_{A[G_v]}(\mathcal O_K/w^e)=(P_v(0))=(Nv).$$
\medskip

Let $K' := K^{I_v}$ denote the maximal sub-extension of $K/F$, which is unramified at $v.$ Let $w'$ denote the prime in $K'$ lying below $w$.
Since $K'/F$ is unramified at $v$, Proposition \ref{prop-unr} leads to an equality of ideals
$${\rm Fitt}_{A[G_v/I_v]}(E(\mathcal{O}_{w'}/w')) = (P_{v}({\sigma}'_v))$$
where ${\sigma}'_v$ denotes the Frobenius element of $v$ in $K'/F$. Since $\mathcal{O}_{w'}/w' \simeq \mathcal{O}_w/w$ (as $K_w/K'_{w'}$ is totally ramified at $w'$), we also have an equality of ideals 
\begin{equation}\label{Fitt-ev}{\rm Fitt}_{A[G_v/I_v]}(E(\mathcal{O}_{w}/w)) = (P_{v}({\sigma}'_v))\end{equation}
For simplicity, let $e = |I_v|.$ Then, Proposition A.5.1. in \cite{FGHP} gives isomorphisms of $A[G_v]$--modules
$$e_v(\mathcal{O}_w/w^e) \simeq \mathcal{O}_w/w \hspace{10pt} \text{ and } \hspace{10pt} (1 - e_v)(\mathcal{O}_w/w^e) \simeq w/w^e.$$
Since $\mathcal O_w/w^e=e_v(\mathcal{O}_w/w^e)\oplus(1-e_v)(\mathcal{O}_w/w^e)$, the isomorphisms above lead to an equality of $A[G_v]$--ideals 
\begin{equation}\label{Fitt-product}
    {\rm Fitt}_{A[G_v]}(E(\mathcal{O}_w/w^e))={\rm Fitt}_{A[G_v]}(E(\mathcal{O}_w/w))\cdot {\rm Fitt}_{A[G_v]}(E(w/w^e)).\end{equation}

\begin{rem}
     The definition of Fitting ideals implies that if $M$ is a finitely generated $A[G_v]$--module, then 
$${\rm Fitt}_{A[G_v]}(e_vM)=e_v{\rm Fitt}_{A[G_v]}(M)\oplus (1-e_v)A[G_v], \qquad {\rm Fitt}_{A[G_v]}((1-e_v)M)=e_vA[G_v]\oplus (1-e_v){\rm Fitt}_{A[G_v]}(M).$$
\end{rem}
\medskip

We will use equality \eqref{Fitt-product} to prove \eqref{Fitt-local}. The obvious ring isomorphism $e_vA[G_v] \simeq A_v[G_v/I_v]$ which sends $e_v\sigma_v\to \sigma_v'$ and equality \eqref{Fitt-ev} give 
\begin{equation}\label{Fitt-e_v-again}
{\rm Fitt}_{A[G_v]}(E(\mathcal{O}_w/w)) = e_vP(\sigma_v)A[G_v]\oplus (1-e_v)A[G_v].\end{equation}
Next, we will calculate  ${\rm Fitt}_{A[G_v]}(E(w/w^e)).$ For that, we need a definition and a couple of lemmas.

\begin{Def} Let $R$ be a commutative ring and $M$ a finitely presented $R$--module. $M$ is called quadratically presented if there exists an $n>0$ and an exact sequence of $R$--modules
$$R^n\to R^n\to M\to 0.$$
$M$ is called locally quadratically presented if $M_{\frak m}$ is a quadratically presented $R_{\frak m}$--module for all $\frak m\in{\rm MSpec}(R)$.
\end{Def}
The following computationally useful result is due to Johnston and Nickel. (See \cite{Johnston-Nickel}).
\begin{Lemma}[Johnston--Nickel]\label{Johnston-Nickel}
If $R$ is a commutative ring, $C$ is a locally quadratically presented $R$ and 
$$0\to A\to B\to C\to 0$$
is an exact sequence of finitely presented $R$--modules, then 
$${\rm Fitt}_R(B)={\rm Fitt}_R(A)\cdot {\rm Fitt}_R (C).$$
\end{Lemma}
A supply of locally quadratically presented $A[G_v]$--modules in the context at hand is given by the following Lemma, whose proof uses concepts of group cohomology discussed in detail in the Appendix of \cite{FGHP}.

\begin{Lemma}\label{ct-Lemma} Under the above assumptions, the following hold.
\begin{enumerate}
 \item The $A[G_v]$--submodule $w^i$ of $\mathcal O_w$ is $G_v$--cohomologically trivial ($G_v$--c.t.), for all $i\geq 0$.
 \item For all $j>i\geq 0$, the $A[G_v]$--modules $w^i/w^j$ and $E(w^i/w^j)$ are locally quadratically presented.
\end{enumerate}
\end{Lemma}
\begin{proof}
(1) Fix $i\geq 0$. Since the $A[G_v]$--module $w^i$ is annihilated by $p$ and $p\nmid |I_v|$, we have
$${\rm H}^r(I_v, w^i)=0, \text{ for all }r\geq 1.$$
Consequently, the inflation--restriction sequences in group--cohomology give us exact sequences of $A$--modules
$$0\to {\rm H}^r(G_v/I_v, (w^i)^{I_v})\to {\rm H}^r(G_v, w^i)\to {\rm H}^r(I_v, w^i)^{G_v},$$
for all $r\geq 1$. However, $(w^i)^{I_v}=(w')^{i'}$, for some $i'\geq 0$. Since $K'_{w'}/F_v$ is unramified, we have isomorphisms of $\mathcal O_v[G_v/I_v]$--modules 
$$\mathcal O_v[G_v/I_v]\simeq \mathcal O_{w'}\overset{\pi_v^{i'}}\longrightarrow {w'}^{i'}.$$
where the second isomorphism is multiplication with $\pi_v^{i'}$, where $\pi_v$ is a generator of $v$ (as an ideal of $\mathcal O_v$). This shows that the $\mathcal O_v[G_v/I_v]$--module $(w^i)^{I_v}$ is $G_v/I_v$--induced and therefore $G_v/I_v$--c.t.
Now, the exact sequence above implies that 
$${\rm H}^r(G_v, w^i)=0, \qquad \text{ for all } r\geq 1,$$
which shows that, indeed, $w^i$ is $G_v$--c.t., for all $i\geq 0$.   

(2) Part (1) implies that the $A[G_v]$ modules $w^1/w^j$ are $G_v$--c.t. (since $w^j$ and $w^i$ are. Therefore, since $A$ is a PID, they have projective dimension at most $1$ over $A[G_v]$. However, $w^i/w^j$ is finite and since $A$ is infinite, $w^i/w^j$ must  have projective dimension exactly $1$ over $A[G_v]$. Therefore, there are $\alpha_{\frak m}, \beta_{\frak m}\in\Bbb Z$, with $0\leq \alpha_{\frak m}\leq \beta_{\frak m}$ and exact sequences of $A[G_v]_{\frak m}$--modules 
$$0\to A[G_v]^{\alpha_{\frak m}}\to A[G_v]^{\beta_{\frak m}}\to (w^i/w^j)_{\frak m}\to 0,$$
for all $\frak m\in{\rm MSpec}(A[G_v])$. Again, by the finiteness of $(w^i/w^j)$, one concludes that $\alpha_{\frak m}=\beta_{\frak m}$, which gives the desired local quadratic presentations, concluding the proof.
\end{proof}

\begin{Cor}\label{Fitt(1-e_v)-Cor} Under the above assumptions, the following hold.
\begin{enumerate}
\item For all $i>0$, the identity map induces an isomorphism of $A[G_v]$--modules
$$w^i/w^{i+1}\simeq E(w^i/w^{i+1}).$$
\item We have equalities of $A[G_v]$--ideals 
$${\rm Fitt}_{A[G_v]} E(w/w^e)={\rm Fitt}_{A[G_v]} (w/w^e)= e_vA[G_v]\oplus (1-e_v)Nv\cdot A[G_v].$$
\end{enumerate}
    
\end{Cor}
\begin{proof}
(1) The identity map induces an isomorphism of $\Bbb F_q[G]$--modules between the two modules in question. However, since $\tau(w^i)\subseteq w^{qi}\subseteq w^{i+1}$, $\tau$ acts trivially on $w^i/w^{i+1}$. Therefore, we have 
$$\phi_E(t)\cdot x= t\cdot x + a_1\tau(x) + ... + a_r\tau^r(x)=t\cdot x,\qquad\text{for all }x\in w^i/w^{i+1},$$
which shows that the identity is also an isomorphism of $\Bbb F_q[G_v][t]=A[G_v]$--modules.
\medskip

(2) Consider the short exact sequences of $A[G_v]$--modules 
$$0\to w^i/w^{i+1}\to w/w^{i+1}\to w/w^i\to 0, \qquad 0\to E(w^i/w^{i+1})\to E(w/w^{i+1})\to E(w/w^i)\to 0, $$
for all $i=1, \dots, e-1.$ Apply Lemma \ref{ct-Lemma}(2) and Lemma \ref{Johnston-Nickel} to conclude that 
$${\rm Fitt}_{A[G_v]}(w/w^{i+1})={\rm Fitt}_{A[G_v]}(w/w^{i})\cdot {\rm Fitt}_{A[G_v]}(w^i/w^{i+1}),$$ $${\rm Fitt}_{A[G_v]}E(w/w^{i+1})={\rm Fitt}_{A[G_v]}E(w/w^{i})\cdot {\rm Fitt}_{A[G_v]}E(w^i/w^{i+1}),
$$
for all $i=1, \dots, e-1.$. Consequently, we obtain the following equalities of $A[G_v]$--ideals.
$${\rm Fitt}_{A[G_v]}(w/w^e)=\prod_{i=1}^{e-1}{\rm Fitt}_{A[G_v]}(w^i/w^{i+1}),\qquad 
{\rm Fitt}_{A[G_v]}E(w/w^e)=\prod_{i=1}^{e-1}{\rm Fitt}_{A[G_v]}E(w^i/w^{i+1}).
$$
Now, we combine the above equalities with the $A[G_v]$--module isomorphisms in part (1) to obtain
$${\rm Fitt}_{A[G_v]}E(w/w^e)= {\rm Fitt}_{A[G_v]}(w/w^e).$$
Now, since $w/w^e=(1-e_v)(\mathcal O_w/w^e)=(1-e_v)(\mathcal O_K/w^e)$, we have 
$${\rm Fitt}_{A[G_v]}E(w/w^e)=e_vA[G_v]\oplus (1-e_v){\rm Fitt}_{A[G_v]}\mathcal O_K/w^e=e_vA[G_v]\oplus (1-e_v)Nv\cdot A[G_v],$$
which concludes the proof of the Corollary.
\end{proof}
Now, we use Corollary \ref{Fitt(1-e_v)-Cor}(2) and equalities \eqref{Fitt-product}--\eqref{Fitt-e_v-again} to conclude that 
$${\rm Fitt}_{A[G_v]}E(\mathcal O_w/w)=e_vP(\sigma_v)\cdot A[G_v]\oplus (1-e_v)Nv\cdot A[G_v].$$ However, note that if $P_v(X)=X^r+a_{r-1}X^{r-1}+\dots +a_1X^r+a_0$, then $P_v(0)=a_0=\rho\cdot Nv$ and therefore 
$$P_v(e_v\sigma_v)=e_vP_v(\sigma_v)+(1-e_v)P_v(0)=e_vP_v(\sigma_v)+(1-e_v)\rho N_v.$$
Now, recalling that $\rho\in\Bbb F_q^\times$, we combine the last two displayed equalities to conclude that 
$${\rm Fitt}_{A[G_v]}E(\mathcal O_w/w)=P_v(e_v\sigma_v)A[G_v].$$
This proves equality \eqref{Fitt-local} and concludes the proof of Proposition \ref{prop-tame}.\end{proof}

\vspace{20pt}
\section{The $t$--motive and local $\Bbb F_q$--shtukas associated to $\overline E$}\label{motive-section}
This section is dedicated to the proof of a more precise version of Proposition \ref{g-divides-f-prop}, which will be stated below. The arguments below ensued from a lengthy and enlightening e-mail exchange \cite{private} between the first author and Urs Hartl, to whom we would like to express our gratitude.

As before, $\overline E$ is the reduction modulo $v$ of the original Drinfeld module $E$ defined on $A$ and over $\mathcal O_F$. So, $\overline E$ is a Drinfeld module of rank $r$, defined on $A$, over the finite field $\mathcal O_F/v=\Bbb F_{q^{d_v}}$, and of characteristic $w_0=v\cap A$. In \S\ref{reduction-section}, we denoted by $h$ the height of $\overline E$, picked a prime $v_0\in{\rm MSpec}(A)\setminus\{w_0\}$,  and considered 
$$f_{\overline E}(X):={\rm det}_{A_{v_0}}(X\cdot I_r-\tau_1\mid T_{v_0}(\overline E)), \qquad g_{\overline E}(X):={\rm det}_{A_{w_0}}(X\cdot I_{r-h}-\tau_1\mid T_{w_0}(\overline E)),$$
where $\tau_1:={\rm Frob}_{q^{d_v}}$ is the $q^{d_v}$--power Frobenius endomorphism of $\overline E$, acting naturally and linearly on its various Tate modules. In the same section, we stated the well known result that $f_{\overline E}(X)$ is independent on $v_0\ne w_0$ and in fact $f_{\overline E}(X)\in A[X]$ (which will also be proved very briefly below) and claimed without a proof that $g_{\overline E}(X)\mid f_{\overline E}(X)$ in $A_{w_0}[X]$. The precise goal of this section is to prove the following.

\begin{proposition}\label{g-divides-f-prop-precise} Let $\overline {A_{w_0}}$ be the absolute integral closure 
of $A_{w_0}$. Then the following hold.

\begin{enumerate}
    \item The roots of $g_{\overline E}(X)=0$ in $\overline {A_{w_0}}$ are exactly those roots $\alpha$ of $f_{\overline E}(X)=0$ which satisfy
$$|\alpha|_{w_0}=1,$$
where $|\cdot |_{w_0}$ is the unique extension to $\overline {A_{w_0}}$ of the normalized absolute value on $A_{w_0}$. 
\item If $\alpha$ is such a root, then its multiplicities in $g_{\overline E}(X)$ and $f_{\overline E}(X)$ are the same.

\item In particular, $g_{\overline E}(X)\mid f_{\overline E}(X)$ in $A_{w_0}[X].$
\end{enumerate}
\end{proposition}

Obviously, part (3) of the above Proposition is a direct consequence of parts (1) and (2). In order to prove parts (1) and (2), we need to introduce first the $t$--motive 
$\overline M:=M_{\overline E}$ associated to $\overline E$. We follow most of the notations and conventions in \cite{Goss}, Chapter 5. Let $L_1:=\mathcal O_F/v$ and let $L:=\overline{L_1}=\overline{\Bbb F_q}$ be a fixed algebraic closure. In what follows,  $\Bbb G_a$ denotes the additive affine line, viewed as a scheme over ${\rm Spec}(\Bbb F_q)$. We think of $\overline E$ as a functor from the category of $L_1$--algebras to the category of $A$--modules
$$\overline E: [\text{$L_1$--alg}]\longrightarrow [\text{$A$--mod}], \qquad L'\to \Bbb G_a(L'),$$
where $\Bbb G_a(L')$ is endowed with a natural $A$--module structure via the $\Bbb F_q$--algebra (injective) morphism
$$A\overset{\phi_{\overline E}}{\longrightarrow}L_1\{\tau\}\subseteq L'\{\tau\}={\rm End}_{\Bbb F_q}^{L'}(\Bbb G_a).$$
\begin{Def} As $L$ is a perfect field containing $L_1$ (the field of definition of $\overline E$), we follow loc.cit. and define the $t$--motive over $L$ associated to $E$ as the left $L\{\tau\}\otimes_{\Bbb F_q}A=L\{\tau\}[t]$--module
$$\overline M(L):={\rm Hom}_{\Bbb F_q}^L(\overline E(L), \Bbb G_a(L))=L\{\tau\},$$
endowed with the left $L\{\tau\}\otimes_{\Bbb F_q}A$--module structure given by
$$(\lambda\otimes a)\ast \mu:=\lambda\circ\mu\circ\phi_{\overline E}(a), \qquad\text{ for all }\lambda\in L\{\tau\}, a\in A, \mu\in \overline M(L).$$
\end{Def}

\begin{rem}\label{properties-remark} It is important to note that the $L\{\tau\}[t]$--module $\overline M(L)$ has some distinctive properties (see loc. cit. for proofs): First, it is obvious that $\overline M(L)$ is a free $L\{\tau\}=(L\{\tau\}\otimes 1)$--module of rank $1$ (which is {\bf the dimension} of the $t$--motive $\overline M(L)$) and (less obvious) that it is a free $L[t]=(L\otimes_{\Bbb F_q}A)$--module of rank $r$ (which is {\bf the rank} of the $t$--motive $\overline M(L)$.) Second, it is important to note that since $L$ is perfect, $\tau \overline M(L)$ is an $L\{\tau\}[t]$--submodule of $\overline M(L)$ and, as a consequence of the definition of $\phi_{\overline E}$, we have
$$(1\otimes t-\iota(t)\otimes 1)(\overline M(L)/\tau\overline M(L))=0,$$
where $i:A\to L_1\subseteq L$ is the obvious $\Bbb F_q$--algebra map of kernel $w_0$. 
\end{rem}
It is not difficult to check that the evaluation (perfect) pairing 
$$\overline E(L)\times \overline M(L)\to \Bbb G_a(L), \qquad (e, \mu)\to \mu(e)$$
gives rise to an isomorphism of $A$--modules
$$\xi: \overline E(L)\simeq {\rm Hom}_{L\{\tau\}[t]}(\overline M(L), L((t^{-1}))/tL[t]), \qquad 
e\to \big [\,  \mu\to \overline{\sum_{i\geq 0}\mu(\phi_{\overline E}(t^i)(e))\cdot t^{-i}}\, \big ],$$
where $\tau$ acts on $L((t^{-1}))/tL[t]$ by raising the coefficients of the Laurent series in question to the $q$--th power and $L[t]$ acts via multiplication. For every $f\in A$, this leads to a natural isomorphism 
of $A/f$--modules
$$\xi[f]: \overline E[f]\simeq {\rm Hom}_{L\{\tau\}\otimes_{\Bbb F_q} A/f}\big(\overline M(L)/f,\,L[t]/fL[t]\big),$$
after identifying $L[t]/fL[t]\simeq (L((t^{-1}))/tL[t])[f]$ via the isomorphism $\widehat \rho\to \widehat{t\rho/f}$.
\medskip

Now, we fix an arbitrary $v_0\in{\rm MSpec}(A)$ and let $\pi_{v_0}\in A$ denote the monic generator of $v_0$. We let 
$$A^{nr}_{v_0}:=L\widehat{\otimes}_{\Bbb F_q}A_{v_0}:={\varprojlim_n}\, (L\otimes_{\Bbb F_q}A/v_0^n), \qquad \overline M(L)_{v_0}:=\overline M(L)\widehat\otimes_{A}A_{v_0}:={\varprojlim_n}\, (\overline M(L)\otimes_{A}A/v_0^n).$$
Note that if $d_{v_0}:=[A/v_0:\Bbb F_q]$, then we have natural isomorphisms of topological rings

$$A_{v_0}\simeq {\Bbb F}_{q^{d_{v_0}}}[[\pi_{v_0}]], \qquad A^{nr}_{v_0}\simeq L[[\pi_{v_0}]]^{d_{v_0}}.$$
Further, note that since $\overline M(L)$ is a free $(L\otimes_{\Bbb F_q}A=L[t])$--module
of rank $r$ (see the Remark above), then  $\overline M(L)_{v_0}$ is a free $A_{v_0}^{nr}$--module of rank $r$ and, consequently, a free $L[[\pi_{v_0}]]$--module of rank $rd_{v_0}$. In addition, if we view ${\rm Frob}_q$ as the canonical topological generator of $Gal(L/\Bbb F_q)=Gal(A_{v_0}^{nr}/A_{v_0})$, then the free $A_{v_0}^{nr}$--module $\overline M(L)_{v_0}$ is endowed with a ${\rm Frob}_q$--semilinear endomorphism, abusively denoted $\tau$, and given by
$$\tau:=\tau\widehat\otimes 1: \overline M(L)\widehat\otimes_{A}A_{v_0}\longrightarrow \overline M(L)\widehat\otimes_{A}A_{v_0}.$$

\begin{Def} The data $(\overline M(L)_{v_0}, \tau)$ consisting of the free $A_{v_0}^{nr}$--module $\overline M(L)_{v_0}$ of rank $r$ together with its ${\rm Frob}_q$--semilinear endomorphism $\tau$ defined above is called the local $\Bbb F_q$--shtuka over $L$ associated to $\overline E$ at $v_0$.
\end{Def}

The link between the local shtuka $(\overline M(L)_{v_0}, \tau)$ and the Tate module $T_{v_0}(\overline E)$ is obtained by taking the projective limit as $n\to\infty$ of the 
isomorphisms $\xi[\pi_{v_0}^n]$ defined above, to get an isomorphism of $A_{v_0}$--modules
$${\xi}_{v_0}^{nr}: T_{v_0}(\overline E)\simeq{\rm Hom}_{A_{v_0}^{nr}\{\tau\}}(\overline M(L)_{v_0}, A_{v_0}^{nr}), \qquad \text{ for all }v_0\in{\rm MSpec}(A).$$
The following useful result is due to Hartl--Singh \cite{Hartl-Singh}, building upon earlier work of Laumon \cite{Laumon}.
\begin{proposition}
For all $v_0\in{\rm MSpec}(A)$, the local $\Bbb F_q$--shtuka $(\overline M(L)_{v_0}, \tau)$ over $L$ splits canonically as a direct sum of local $\Bbb F_q$--shtukas over $L$
$$(\overline M(L)_{v_0},\tau)=(\overline M(L)_{v_0}^{et}, \tau)\oplus (\overline M(L)_{v_0}^{nil}, \tau),$$
where  $\overline M(L)_{v_0}^{et}$ is the maximal $A_{v_0}^{nr}\{\tau\}$--submodule of $\overline M(L)_{v_0}$ on which the restriction of $\tau$ is bijective and $\overline M(L)_{v_0}^{nil}$ is the maximal $A_{v_0}^{nr}\{\tau\}$--submodule of $\overline M(L)_{v_0}$ on which the restriction of $\tau$ is topologically nilpotent (i.e. there exists an $n>0$ such that $\tau^n (\overline M(L)^{nil})\subseteq \pi_{v_0}\overline M(L)^{nil}$.)
\end{proposition}
\begin{proof} See the proofs of Propositions 2.7 and 2.9 in \cite{Hartl-Singh} and keep in mind that the field $L$ is perfect in our context.
Also, note that in loc.cit. the authors work with the $L[[\pi_{v_0}]]$--module structures (local) rather than the $A_{v_0}^{nr}$---module structures (semi-local), but the transition local/semi--local is seamless. 
\end{proof}
In light of the above result, since obviously $\tau$ acts as an isomorphism of $A_{v_0}^{nr}$ (i.e. $(A_{v_0}^{nr}, \tau)$ is an \'etale local $\Bbb F_q$--shtuka over $L$), the isomorphism of $A_{v_0}$--modules $\xi_{v_0}^{nr}$ defined above can be rewritten as 
$${\xi}_{v_0}^{nr}: T_{v_0}(\overline E)\simeq{\rm Hom}_{A_{v_0}^{nr}\{\tau\}}(\overline M(L)^{et}_{v_0}, A_{v_0}^{nr}), \qquad \text{ for all }v_0\in{\rm MSpec}(A).$$
The next crucial step is provided by the following extension to the case of ${\rm GL_n}(A_{v_0}^{nr})$ of Lang's well known theorem 
on ${\rm GL}_n(L)$ (see \cite{Lang}).

\begin{Lemma}[Hartl, \cite{private}]\label{gen-Lang-lemma} Under the above assumptions, the following hold.
\begin{enumerate}
\item The map ${\rm GL}_n(A_{v_0}^{nr})\to {\rm GL_n}(A_{v_0}^{nr})$ taking $X\to {\rm Frob}_q(X)^{-1}\cdot X$ is surjective.
\item Any free $A_{v_0}^{nr}$--module $\mathcal M$ of finite rank $n$, endowed with a bijective, ${\rm Frob}_q$--semilinear endomorphism $t$ satisfies the property that the standard map
$$\mathcal M^{t=1}\otimes_{A_{v_0}}A_{v_0}^{nr}\to \mathcal M, \qquad a\otimes m\to am$$
is an isomorphism of $A_{v_0}^{nr}$--modules.
\end{enumerate}
\end{Lemma}
\begin{proof}
Since we have a ring isomorphism $A_{v_0}^{nr}\simeq L[[\pi_{v_0}]]^{d_{v_0}}$, it suffices to prove part (1) for ${\rm GL}_n(L[[\pi_{v_0}]])$. So, given a matrix $A\in{\rm GL}_n(L[[\pi_{v_0}]])$, i.e. 
$$A=A_0+A_1\cdot\pi_{v_0}+\dots , \qquad \text{ with } A_0\in{\rm GL}_n(L)\text{ and } A_i\in M_n(L),\text{ for all } i\geq 1,$$
we need to find a matrix $X\in {\rm GL}_n(L[[\pi_{v_0}]]$ given by 
$$X=X_0+X_1\cdot\pi_{v_0}+\dots , \qquad \text{ with } X_0\in{\rm GL}_n(L)\text{ and } X_i\in M_n(L),\text{ for all } i\geq 1,$$
such that the matrices $X_i$ satisfy the relations
$$\sum_{i=0}^m{\rm Frob}_q(X_{m-i})\cdot A_i=X_m,\qquad \text{ for all }m\geq 0.$$
Lang's theorem (see loc.cit.) implies that part (1) is true for ${\rm GL}_n(L)$, so we can find a matrix $X_0\in{\rm GL}_n(L)$ satisfying the $0$--th relation above. After multiplying the $m$--th relation above to the right by $A_0^{-1}{\rm Frob}_q(X_0)^{-1}=X_0^{-1}$ we obtain the equivalent relation
$$X_mX_0^{-1} = {\rm Frob}_q(X_mX_0^{-1})+ \sum_{i\geq 1}^m{\rm Frob}_q(X_{m-i})\cdot A_i\cdot X_0^{-1},$$
which consists of one Artin--Schreier equation for each entry of $X_mX_0^{-1}$. Since $L$ is algebraically closed, these equations have solutions. Therefore, inductively, one can find matrices $X_m$, for all $m\geq 0$, as desired.

Part (2) follows immediately from part (1) in a standard way: take a basis $\overline e$
of $\mathcal M$ over $A_{v_0}^{nr}$ and let $A$ be the matrix of $t$ in that basis. Let $X\in {\rm GL}_n(A_{v_0}^{nr})$ such that $A={\rm Frob}_q(X)^{-1}\cdot X$. Then $\overline {e'}:=X\cdot\overline{e}$ is an $A_{v_0}^{nr}$-basis of $\mathcal M$ which is contained in $\mathcal M^{t=1}$. This concludes the proof.    
\end{proof}

By applying the Lemma above to $\mathcal M:=\overline M(L)_{v_0}^{et}$ and $t=\tau$, we conclude that we have the following natural isomorphisms of $A_{v_0}^{nr}\{\tau\}$--modules
$$\overline M(L)_{v_0}^{et}\simeq (\overline M(L)_{v_0}^{et})^{\tau=1}\otimes_{A_{v_0}}A_{v_0}^{nr}, \qquad \text{ for all }v_0\in{\rm MSpec}(A).$$
The above isomorphism leads to a further isomorphism of $A_{v_0}$--modules
$${\rm Hom}_{A_{v_0}^{nr}\{\tau\}}(\overline M(L)^{et}_{v_0}, A_{v_0}^{nr})\simeq {\rm Hom}_{A_{v_0}}((\overline M(L)^{et}_{v_0})^{\tau=1}, A_{v_0}), \qquad \text{ for all }v_0\in{\rm MSpec}(A),$$
which, if composed with the map $\xi_{v_0}^{nr}$ gives an isomorphism of $A_{v_0}$--modules
$$\xi_{v_0}: T_{v_0}(\overline E)\simeq {\rm Hom}_{A_{v_0}}((\overline M(L)^{et}_{v_0})^{\tau=1}, A_{v_0}), \qquad \text{ for all }v_0\in{\rm MSpec}(A).$$
This prompts the following.
\begin{Def}\label{cohomology}
    The first $v_0$--adic \'etale cohomolgy group of $\overline E$ is defined by
$${\rm H}^1_{et}(\overline E, \, A_{v_0}):= (\overline M(L)^{et}_{v_0})^{\tau=1}, \qquad \text{ for all }v_0\in{\rm MSpec}(A).$$
{\bf Note} that the maps $\xi_{v_0}$ lead to the following $A_{v_0}$--module isomorphisms.
$$T_{v_0}(\overline E)^\ast:= {\rm Hom}_{A_{v_0}}(T_{v_0}(E), A_{v_0})\simeq {\rm H}^1_{et}(\overline E, \, A_{v_0}), \qquad\text{ for all }v_0\in{\rm MSpec}(A).$$
Further, the first $v_0$--adic crystalline cohomology group of $\overline E$ is defined by
$${\rm H}^1_{cris}(\overline E, A_{v_0}^{nr}):=\overline M(L)_{v_0}.$$
{\bf Note} that for all $v_0\in{\rm MSpec}(A)$ we have isomorphisms and inclusions of $A_{v_0}^{nr}$--modules
$$T_{v_0}(\overline E)^\ast\otimes_{A_{v_0}}A_{v_0}^{nr} ={\rm H}^1_{et}(\overline E, A_{v_0})\otimes_{A_{v_0}}A_{v_0}^{nr}\simeq {\rm H}^1_{cris}(\overline E, A_{v_0}^{nr})^{et}\subseteq {\rm H}^1_{cris}(\overline E, A_{v_0}^{nr}).$$
\end{Def}
\medskip

The following holds at primes $v_0$ different from the characteristic of $\overline E$. (See \cite{Goss}, Chapter 5 as well.)
\begin{Lemma}\label{etale-non-char-lemma}
 If $v_0$ is an element in ${\rm MSpec}(A)$ different from the characteristic $w_0$ of $\overline E$, then
 $$\overline M(L)_{v_0}^{et}=\overline M(L)_{v_0},\qquad (\overline M(L)_{v_0})^{\tau=1}\otimes_{A_{v_0}}A_{v_0}^{nr}\simeq \overline M(L)_{v_0}.$$
 In other words, $\tau$ is bijective on $\overline M(L)_{v_0}$ and we have canonical isomorphisms of $A_{v_0}^{nr}[\tau_1]$--modules
 $$T_{v_0}(\overline E)^\ast\otimes_{A_{v_0}}A_{v_0}^{nr}\simeq {\rm H}^1_{et}(\overline E, A_{v_0})\otimes_{A_{v_0}}A_{v_0}^{nr}\simeq {\rm H}^1_{cris}(\overline E, A_{v_0}^{nr}).$$
\end{Lemma}
\begin{proof}(sketch) It is easy to show that since $\tau\nmid\phi_{\overline E}(\pi_{v_0}^n)$ in $L\{\tau\}$, $\tau$ is injective and therefore bijective on the finite dimensional $L$--vector spaces $\overline M(L)/\pi_{v_0}^n$, for all $n\geq 1$.
The bijection of $\tau$ on $\overline M(L)_{v_0}$ is obtained now by taking the projective limit as $n\to\infty$.
\end{proof}

Now, we have all the necessary ingredients to prove parts (1) and (2) of Proposition \ref{g-divides-f-prop-precise}.
\begin{proof} First, take $v_0\in{\rm MSpec}(A)\setminus\{w_0\}$ and observe that, based on the previous Lemma and Definition, we have canonical isomorphisms of $A_{v_0}^{nr}[\tau_1]$--modules
$$T_{v_0}(\overline E)^\ast\otimes_{A_{v_0}}A_{v_0}^{nr}\simeq {\rm H}^{1}_{cris}(\overline E, A_{v_0}^{nr})= \overline M(L)\widehat{\otimes}_{A}A_{v_0}\simeq \overline M(L)\otimes_{L\otimes_{\Bbb F_q}A}(L\widehat{\otimes}_{\Bbb F_q}A_{v_0})\simeq \overline M(L)\otimes_{L[t]}A_{v_0}^{nr}.$$
As a consequence, from the definition of $f_{\overline E}$, we have 
$$f_{\overline E}(X)={\rm det}_{A_{v_0}}(XI_r-\tau_1\mid T_{v_0}(\overline E)^\ast)={\rm det}_{A_{v_0}^{nr}}(XI_r-\tau_1\mid {\rm H}^{1}_{cris}(\overline E, A_{v_0}^{nr}))={\rm det}_{L[t]}(XI_r-\tau_1\mid \overline M(L)).$$
The last equality proves that $f_{\overline E}(X)$ is independent of $v_0$ and that it has coefficients in $L[t]$. Further, if one applies  the analogues of Lemmas \ref{gen-Lang-lemma} and \ref{etale-non-char-lemma} to the finite, \'etale $\Bbb F_q$--shtukas $\overline M(L)/{v_0}^n$ over $L$ (see \cite{Hartl-Singh}), one concludes that $f_{\overline E}(X)$ has coefficients in $A_{v_0}$. Since $L[t]\cap A_{v_0}=A$ (intersection viewed inside $A_{v_0}^{nr}$), $f_{\overline E}(X)$ has coefficients in $A$, as stated before.
\medskip

Now, from the definitions, we also have similar natural isomorphisms of $A_{w_0}^{nr}[\tau_1]$--modules
$${\rm H}^{1}_{cris}(\overline E, A_{w_0}^{nr})= \overline M(L)\widehat{\otimes}_{A}A_{w_0}\simeq \overline M(L)\otimes_{L\otimes_{\Bbb F_q}A}(L\widehat{\otimes}_{\Bbb F_q}A_{w_0})\simeq \overline M(L)\otimes_{L[t]}A_{w_0}^{nr}.$$
Therefore, when combining these with the second note in Definition \ref{cohomology}, we obtain equalities
$$\begin{aligned}
    f_{\overline E}(X)&={\rm det}_{L[t]}(XI_r-\tau_1\mid \overline M(L))={\rm det}_{A_{w_0}^{nr}}(XI_r-\tau_1\mid \overline M(L)\otimes_{L[t]}A_{w_0}^{nr})\\
    &=
{\rm det}_{A_{w_0}^{nr}}(XI_r-\tau_1\mid {\rm H}^1_{cris}(\overline E, A_{w_0}^{nr}))\\
&=
{\rm det}_{A_{w_0}^{nr}}(XI_r-\tau_1\mid {\rm H}^1_{cris}(\overline E, A_{w_0}^{nr})^{et})\cdot {\rm det}_{A_{w_0}^{nr}}(XI_r-\tau_1\mid {\rm H}^1_{cris}(\overline E, A_{w_0}^{nr})^{nil})\\
&={\rm det}_{A_{w_0}}(XI_{r-h}-\tau_1\mid T_{w_0}(\overline E)^\ast)\cdot {\rm det}_{A_{w_0}^{nr}}(XI_r-\tau_1\mid {\rm H}^1_{cris}(\overline E, A_{w_0}^{nr})^{nil})\\
&=g_{\overline E}(X)\cdot {\rm det}_{A_{w_0}^{nr}}(XI_r-\tau_1\mid {\rm H}^1_{cris}(\overline E, A_{w_0}^{nr})^{nil}).
\end{aligned}
 $$ 
 Firstly, this shows that $g_{\overline E}(X)$ divides $f_{\overline E}(X)$ in $A_{w_0}^{nr}[X]$. However, since $g_{\overline E}(X), f_{\overline E}(X)$ are both in the polynomial ring $A_{w_0}[X]=A_{w_0}^{nr}[X]^{\tau=1}$, this divisibility holds in $A_{w_0}[X]$. Secondly, note that 
 $$g_{\overline E}(X)={\rm det}_{A_{w_0}^{nr}}(XI_r-\tau_1\mid {\rm H}^1_{cris}(\overline E, A_{w_0}^{nr})^{et}), \qquad f_{\overline E}(X)/g_{\overline E}(X)={\rm det}_{A_{w_0}^{nr}}(XI_r-\tau_1\mid {\rm H}^1_{cris}(\overline E, A_{w_0}^{nr})^{nil}).$$
 Since $\tau_1$ (which is a power of $\tau$) acts bijectively on the \'etale piece and topologically nilpotently on the nilpotent piece of ${\rm H}^1_{cris}(\overline E, A_{w_0}^{nr})$, the roots of $g_{\overline E}$ are all $w_0$--adic units and the roots of $f_{\overline E}(X)/g_{\overline E}(X)$ are in the maximal ideal of an absolute integral closure $\overline {A_{w_0}}$ of $A_{w_0}$. This concludes the proof. 
\end{proof}

\section{Final thoughts: the more general case of pure, abelian $t$-modules}\label{tmodule-section}

Most of the techniques developed in this paper can be extended to the more general case of pure, abelian $t$--modules, as indicated briefly below. The very general case of abelian $t$--modules, whose equivariant $L$--functions are the main object of study in \cite{Green-Popescu}, will be treated in upcoming work.
\medskip

With notations as in \S\ref{statement-section}, let $M_n(\mathcal{O}_F)$ be the ring of $n \times n$ matrices with entries in $\mathcal{O}_F$, for a fixed $n\geq 1$. 

\begin{Def} A $t$--module  $E$ of dimension $n$, defined over $\mathcal O_F$ is an $\mathbb{F}_q$-algebra morphism 
$$\phi_E: A \to M_n(\mathcal{O}_F\{\tau\}), \qquad\phi_E(t) = M_0\tau^0 + M_1\tau + ... + M_\ell \tau^\ell,$$
where $M_i \in M_n(\mathcal{O}_F)$, $\ell\geq 1$, $M_{\ell}\ne 0$ and $(M_0 - t\cdot I_n)^n=0$. 
\end{Def}
\medskip

For any $\mathcal O_F\{\tau\}$--module $X$, one denotes by 
$$E(X), \qquad {\rm Lie}_E(X)$$ the additive group $X^{\oplus n}$, endowed with the $A$--module structure given by $\phi_E$ and 
${\rm ev}_{\tau=0}\circ \phi_E$ (evaluation of $\phi_E(a)$ at $\tau=0$, for all $a\in A$), respectively.
\medskip

Note that Drinfeld modules defined on $A$ and over $\mathcal O_F$ are the $1$--dimensional $t$--modules defined in \S\ref{statement-section} above. However, unlike the $1$--dimensional case, where the rank of the Drinfeld module is simply equal to $\ell$, in the case of higher dimensional $t$--modules, introducing the notion of rank is much more subtle. We indicate briefly how this is done below. This will lead to the more restrictive notion of an abelian $t$--module.

Let ${\overline F}$ be the separable closure of $F$. As in \S\ref{motive-section}, we let $E(\overline F):=\mathbb G_a^n(\overline F)$ (where $\mathbb G_a$ is viewed as a group scheme over $\Bbb F_q$), endowed with the $A$--module structure induced by the structural morphism
$$\phi_E: A \to M_n(\mathcal{O}_F)\{\tau\}\subseteq M_n(\overline F)\{\tau\}={\rm End}^{\overline F}_{\Bbb F_q}(E(\overline F)).$$
As in \S\ref{motive-section}, we define the $t$--motive over $\overline F$ associated to $E$ by
$$M(\overline F):=M_E(\overline F):={\rm Hom}_{\Bbb F_q}^{\overline F}(E(\overline F), \Bbb G_a(\overline F)),$$
endowed with the $(\overline F\{\tau\}\otimes_{\Bbb F_q}A)=\overline F\{\tau\}[t]$--module structure given by
$$(\lambda\otimes a)\ast\mu:=(\lambda\circ\mu\circ\phi_{E}(a)), \qquad \text{ for all } a\in A, \lambda\in \overline F\{\tau\}, \mu\in M_E(\overline F).$$
As an $\overline F\{\tau\}$--module, $M_E(\overline F)$ is clearly free of rank $n$ (which is called the dimension of the motive). However, as an $L[t]$--module this is not finitely generated, in general, which brings us to the following.
\begin{Def}
The $t$--module $E$ and the $t$--motive $M_E(\overline F)$ are called abelian if $M_E(\overline F)$ is a finitely generated and (automotically, see \cite{Goss} Chapter 5) free $\overline F[t]$--module. In that case, the rank over $\overline F[t]$ of $M_E(\overline F)$ is called the rank of $E$ and also called the rank of $M_E(\overline F)$.   
\end{Def}
It is not hard to see that if $E$ is a Drinfeld module of classical rank $r$, then $M_E(\overline F)$ is $\overline F[t]$--free of rank $r$. So, Drinfeld modules are $1$--dimensional abelian $t$--modules and the two notions of rank coincide.

Assuming now that $E$ is an abelian $t$--module as above, of dimension $n$ and rank $r$, a prime $v\in{\rm MSpec}(\mathcal O_F)$ is called a prime of good reduction for $E$ if the $t$--module $\overline E$ of dimension $n$, defined over $\mathcal O_F/v$ by
$$\mathcal \phi_{\overline E}: A\overset{\phi_E}\longrightarrow M_n(\mathcal O_F\{\tau\})\twoheadrightarrow M_n(\mathcal O_F/v\{\tau\})$$
(composing $\phi_E$ with the usual mod $v$ reduction map) is abelian of the same rank $r$ as $E$.
    
Now, assume that $E$ is an abelian $t$--module as above, of dimension $n$ and rank $r$, and let $v_0\in {\rm MSpec}(A)$ and $v\in{\rm MSpec}(\mathcal O_F)$, with $v\nmid v_0$ and such that $E$ has good reduction at $v$. One can show without difficulty (see \cite{Goss} Chapter 5) that one has isomorphisms of $A_{v_0}$--modules, for all $v_0\in{\rm MSpec}(A).$
$$E[v_0^m]\simeq (A_{v_0}/v_0^m)^r, \qquad T_{v_0}(E):=\varprojlim_m E[v_0^m]\simeq A_{v_0}^r.$$
Also, the $G(\overline F/F)$--representation $T_{v_0}(E)$ is unramified at all $v\in\text{MSpec}(\mathcal O_F)$ with $v\nmid v_0$ and the polynomials 
$$P_v(X):={\rm det}_{A_{v_0}}(X\cdot I_r-\sigma(v)\mid T_{v_0}(E))$$
are independent of $v$ and have coefficients in $A$. (See \cite{Green-Popescu} and the references therein.) It is also true (see \cite{Green-Popescu}) that if $v$ is tamely ramified in $K/F$, then the finite $A[G]$--modules 
$E(\mathcal O_K/v)$ and ${\rm Lie}_E(\mathcal O_K/v)$ are free over $\Bbb F_q[G]$ and therefore the monic polynomials 
$$|E(\mathcal O_K/v)|_G, \qquad |{\rm Lie}_E(\mathcal O_K/v)|_G$$
are well defined in $A[G]=\Bbb F_q[G][t]$. The abelian $t$--module analogue of Theorem \ref{main-theorem} is the equality
$$\frac{P_v(\sigma_ve_v)}{P_v(0)}=\frac{|E(\mathcal O_K/v)|_G}{|{\rm Lie}_E(\mathcal O_K/v)|_G}.$$
The above equality was used in \cite{Green-Popescu} to prove an equivariant Tamagawa number formula for the value at $0$ of the $G$--equivariant $L$--function associated to the data $(E, F/K)$ as above.
\medskip

In order to prove the above displayed equality in the general case of abelian $t$--modules $E$, one can try to extend to the higher dimensional context the techniques developed in the previous sections of this paper. And, indeed, it is not very difficult to see that the techniques developed in \S\ref{Fitting-Tate-section}, \S\ref{tame-section} and \S\ref{motive-section} extend in full to the higher dimensional context. In particular, Propositions \ref{g-divides-f-prop} and \ref{g-divides-f-prop-precise} hold true in general. However, a higher dimensional analogue of Proposition \ref{eigenvalues-prop} does not hold true, in general. A slightly weaker version of Proposition \ref{eigenvalues-prop} (see Theorem 5.6.10 in \cite{Goss}), which is strong enough for our purposes, holds true for a smaller class of abelian $t$--modules $E$ and primes of good reduction $v\in{\rm MSpec}(\mathcal O_F)$, namely those for which $\overline E$ (the reduction of $E$ mod $v$) is, in addition, a pure $t$--module. The notion of purity is technical and will not be described here. Instead, we refer the interested reader to \S5.5 in \cite{Goss}. To summarize, if we assume purity of $\overline E$, then the last displayed equality above can be proved with the techniques developed in this paper. If, for example, $E$ is a tensor product of Drinfeld modules (in the category of $t$--modules), then purity is achieved. The more general case of arbitrary abelian $t$--motives will be treated in an upcoming paper.

\bibliographystyle{amsplain}

\bibliography{ref}

\providecommand{\bysame}{\leavevmode\hbox to3em{\hrulefill}\thinspace}
\providecommand{\MR}{\relax\ifhmode\unskip\space\fi MR }
\providecommand{\MRhref}[2]{%
  \href{http://www.ams.org/mathscinet-getitem?mr=#1}{#2}
}
\providecommand{\href}[2]{#2}
\begin{thebibliography}{10}

\bibitem{Beaumont}
T.~Beaumont, \emph{On equivariant class number formulas for $t$--modules}, https://arxiv.org/abs/2110.15696, 36 pages.

\bibitem{FGHP}
J.~Ferrara, N.~Green, Z.~Higgins, and C.D. Popescu, \emph{An equivariant {T}amagawa number formula for {D}rinfeld modules and applications}, J. Algebra and Number Theory \textbf{16} (2022), no.~9, 2215--2264.

\bibitem{Goss}
D.~Goss, \emph{Basic structures of function field arithmetic}, Springer Berlin, Heidelberg, 1998.

\bibitem{Green-Popescu}
N.~Green and C.D. Popescu, \emph{An equivariant {T}amagawa number formula for {$t$}--modules and applications}, https://arxiv.org/abs/2206.03541 (June 2022), 27 pages.

\bibitem{GP}
C.~Greither and C.D. Popescu, \emph{The {G}alois module structure of $\ell$-adic realizations of {P}icard 1-motives and applications}, Int. Math. Res. Not. \textbf{2012} (2012), no.~5, 986--1036.

\bibitem{private}
U.~Hartl, \emph{Private email communication with {C}.{D}. {P}opescu},  (2024).

\bibitem{Hartl-Singh}
U.~Hartl and R.K. Singh, \emph{Local shtukas and divisible local {A}nderson modules}, Canad. J. Math. \textbf{71} (2019), no.~5, 1163--1207. \MR{4010425}

\bibitem{Johnston-Nickel}
H.~Johnston and A.~Nickel, \emph{Noncommutative {F}itting invariants and improved annihilation results}, J. Lond. Math. Soc. (2) \textbf{88} (2013), no.~1, 137--160. \MR{3092262}

\bibitem{Lang}
S.~Lang, \emph{Algebraic groups over finite fields}, Amer. J. Math. \textbf{78} (1956), 555--563. \MR{86367}

\bibitem{Laumon}
G.~Laumon, \emph{Cohomology of {D}rinfeld modular varieties. {P}art {I}}, Cambridge Studies in Advanced Mathematics, vol.~41, Cambridge University Press, Cambridge, 1996, Geometry, counting of points and local harmonic analysis. \MR{1381898}

\bibitem{Matsumura}
H.~Matsumura, \emph{Commutative algebra}, The Benjamin/Cummings Publishing Co. Inc., Second Edition, 1980.

\bibitem{Ramachandran-thesis}
N.~Ramachandran, \emph{Some {F}itting ideal computations in {I}wasawa theory over {$\Bbb Q$} and in the theory of {D}rinfeld modules}, Ph.D. thesis, University of California San Diego (2024).

\bibitem{Taelman}
L.~Taelman, \emph{Special {$L$}-values of {D}rinfeld modules}, Ann. of Math. (2) \textbf{175} (2012), no.~1, 369--391. \MR{2874646}

\end{thebibliography}

\end{document}